\documentclass[reqno]{amsart}
\usepackage{graphicx}
\usepackage{amssymb, amsmath, amsthm}
\usepackage{amsfonts}
\usepackage{amssymb}
\usepackage{float}
\usepackage{mathrsfs}
\usepackage{enumitem} 
\newcommand{\dis}{\displaystyle}

\newcommand{\divv}{\text{\rm div}}

\newcommand{\R}{\mathbb R}

\numberwithin{equation}{section}
\numberwithin{theorem}{section}

\numberwithin{figure}{section}

\begin{document}

%
%
%
%
%
%
%
%
%




\title[Blow-up criterion for MHD equations]
 {A blow-up criterion for the 3D compressible magnetohydrodynamics in terms of density}
 
\author{Anthony Suen} 

\address{Department of Mathematics\\Indiana University Bloomington, IN 47405}

\email{cksuen@indiana.edu}


\date{June 17, 2012}



\begin{abstract}

We study an initial boundary value problem for the 3D magnetohydrodynamics (MHD) equations of compressible fluids in $\R^3$. We establish a blow-up criterion for the local strong solutions in terms of the density and magnetic field. Namely, if the density is away from vacuum ($\rho= 0$) and the concentration of mass ($\rho=\infty$) and if the magnetic field is bounded above in terms of $L^\infty$-norm, then a local strong solution can be continued globally in time.

\end{abstract}

\maketitle
\section{Introduction}
We prove a blow-up criteria for the smooth solutions to the compressible magnetohydrodynamics (MHD) in three space dimensions (see Cabannes \cite{cabannes} for a more comprehensive discussion on the system):
\begin{align}
 \rho_t + \divv (\rho u) &=0, \label{1.1}\\
(\rho u^j)_t + \divv (\rho u^j u) + P(\rho)_{x_j} + ({\textstyle\frac{1}{2}}|B|^2)_{x_j}-\divv (B^j B)&= \mu \Delta u^j +
\lambda \, \divv \,u_{x_j},\label{1.2}\\ 
B^{j}_{t} + \divv (B^j u - u^j B) &=\nu\Delta B^j,\label{1.3}\\
\divv\,B &= 0.\label{1.4}
\end{align}
Here $u=(u^1,u^2,u^3)$ and $B=(B^1,B^2,B^3)$ are functions of $x\in\R^3$ and $t\ge0$ representing density, velocity and magnetic field; $P=P(\rho)$
is the pressure; $\varepsilon,\lambda,\nu$ are viscous constants. The system \eqref{1.1}-\eqref{1.4} is solved subjected to some given initial data:
\begin{align}
(\rho,u,B)(x,0)&=(\rho_0,u_0,B_0)(x).\label{1.5}
\end{align}

The local existence of smooth solutions to the MHD system \eqref{1.1}-\eqref{1.4} as well as the global existence of smooth solutions and weak solutions are studied by many mathematicians in decades, see \cite{hu1}, \cite{hu2}, \cite{sart}, \cite{kawashima}, \cite{suenhoff}. When the initial data is taken to be close to a constant state in $H^3(\R^3)$, Kawashima \cite{kawashima} constructed global-in-time $H^3(\R^3)$-solutions. Later, Suen and Hoff \cite{suenhoff} generalized Kawashima's results to obtain global smooth solutions when the initial data is taken to be $H^3(\R^3)$ but only close to a constant state in $L^2(\R^3)$. The existence of global weak solutions to \eqref{1.1}-\eqref{1.4} with large initial data was proved by Hu and Wang \cite{hu1}-\cite{hu2} and Sart \cite{sart} which are extensions of Lions-type weak solutions  \cite{lions} for the Navier-Stokes system. With initial $L^2$-data close to a constant state, Suen and Hoff \cite{suenhoff} generalized Hoff-type intermediate weak solutions \cite{hoff1}-\cite{hoff3} to obtain global solutions to the system \eqref{1.1}-\eqref{1.4}.

On the other hand, the global existence of smooth solution to the MHD system \eqref{1.1}-\eqref{1.4} with arbitrary smooth data is still unknown. For the corresponding Navier-Stokes system, Z. Xin \cite{xin} proved that smooth solution will blow up in finite time in the whole space when the initial density has compact support, while Rozanova \cite{rozanova} showed similar results for rapidly decreasing initial density. Recently Fan-Jiang-Ou \cite{fan} established some blow-up criteria for the classical solutions to 3D compressible flows, which were further extended by Lu-Du-Yao \cite{lu} for MHD system.

The main goal of the present paper is to generalize the corresponding results of Sun-Wang-Zhang \cite{sun} to the MHD system \eqref{1.1}-\eqref{1.4}. When the initial vacuum is allowed, Y. Sun, C. Wang and Z. Zhang obtained a blow-up criterion in terms of the upper bound of the density for the strong solution to the 3-D compressible Navier-Stokes equations. With the presence of magnetic field, we are able to obtain parallel results as in \cite{sun} except that we {\em do not} allow vacuum in the initial density.

\medskip

We now give a precise formulation of our results. First concerning the assumptions on the parameters, we have
\begin{enumerate}[label={\upshape{(1.\arabic*)}}, ref={\upshape{2.}\arabic*}, topsep=0.3cm, itemsep=0.15cm] 
\setcounter{enumi}{5}
\item There exists $K>0$ such that $P(\rho)=K\rho$ for all $\rho>0$;
\end{enumerate}
\setcounter{equation}{6}
\begin{align}\label{1.7}
\varepsilon,\lambda,\nu>0 \text{ and } \lambda<\varepsilon.
\end{align}

\medskip

For the initial data, we assume that
\begin{align}\label{1.8}
\rho_0,u_0,B_0\in H^3(\R^3) \text{ with } \inf(\rho_0)>0 \text{ and } \divv(B_0)=0
\end{align}

and we also write
\begin{equation}\label{9}
C_{0}=\|\rho_0-\tilde\rho\|^2_{H^3}+\|u_0\|^2_{H^3}+\|B_0\|^2_{H^3}.
\end{equation}

\medskip

We make use of the following standard facts (see Ziemer \cite{ziemer} Theorem~2.1.4, Remark~2.4.3, and Theorem~2.4.4, for example). First, given $r\in[2,6]$ there is a constant $C(r)$ such that for $w\in H^1 (\R^3)$,
\begin{equation}\label{1.9}
\|w\|_{L^r(\R^3)} \le C(r) \left(\|w\|_{L^2(\R^3)}^{(6-r)/2r}\|\nabla w\|_{L^2(\R^3)}^{(3r-6)/2r}\right).
\end{equation}

\medskip

We denote the material derivative of a given function $v$ by
$\dot{v}=v_t + \nabla v\cdot u$,
and if $X$ is a Banach space we will abbreviate $X^3$ by $X$. Finally if $I\subset [0,\infty)$ is an interval, $C^1(I;X)$ will be the elements $v\in C(I;X)$ such that the distribution derivative $v_t\in {\mathcal D}'(\R^3\times{\rm int}\,I)$ is realized as an element of $C(I;X)$.
\medskip

We recall a local existence theorem for \eqref{1.1}-\eqref{1.4} by Kawashima \cite{kawashima}, pg.~34--35 and pg.~52--53:

\medskip

\noindent{\bf Theorem~1.1 (Kawashima)} \em Assume that $\varepsilon,\lambda,\nu$ are strictly positive and that the pressure $P$ satisfies {\rm (1.6)}. Then given $\tilde\rho>0$ and $C_3>0$, there is a positive time $T$ depending on $\tilde\rho$, $C_3$ and the parameters $\varepsilon,\lambda,\nu,P$ such that if the initial data $(\rho_0-\tilde\rho,u_0,B_0)$ is given satisfying \eqref{1.8} and 
\begin{align*}
C_0<C_3,
\end{align*}
then there is a solution $(\rho-\tilde,u,B)$ to \eqref{1.1}-\eqref{1.4} defined on $\R^3\times[0,T]$ satisfying
\begin{equation}\label{1.1a}
\rho - \tilde{\rho}\in C([0,T];H^{3}(\R^3))\cap C^1 ([0,T];H^{2}(\R^3))
\end{equation}
and
\begin{equation}\label{1.1b}
u,B\in C([0,T];H^{3}(\R^3))\cap C^1 ([0,T];H^1 (\R^3))\cap L^2([0,T];H^4 (\R^3)).
\end{equation}
\rm
\bigskip
 
The following is the main result of this paper:

\medskip

\noindent{\bf Theorem~1.2} \em Assume that the system parameters satisfy {\rm (1.6)}-\eqref{1.7}. Given $\tilde\rho>0$, suppose $(\rho_0-\tilde\rho,u_0,B_0)$ satisfies \eqref{1.8}. Assume that $(\rho-\tilde\rho,u,B)$ is the smooth solution as constructed in {\rm Theorem~1.1}, and let $T^*\ge T$ be maximal existence time of the solution. If $T^* <\infty$, then we have
\begin{align*}
\lim_{t\rightarrow T^*}||\rho||_{L^\infty((0,t)\times\R^3)}+||\rho^{-1}||_{L^\infty((0,t)\times\R^3)}+||B||_{L^\infty((0,t)\times\R^3)}=+\infty.
\end{align*}
\rm
\bigskip

The rest of the paper is organized as follows. We begin the proofs of Theorem~1.2 in section~2 with a number of {\em a priori} bounds for local-in-time smooth solutions. We make an important use of estimates on the Lam\'{e} operator $L$ which are mainly inspired by \cite{hoff1} and \cite{sun}. Finally in section~3 we prove Theorem~1.2 via a contradiction argument by deriving higher order $H^3$-bounds for smooth solutions. 

\bigskip

\section{A prior estimates}

\bigskip

In this section we derive {\em a prior} estimates for the local solution $(\rho-\tilde\rho,u,B)$ on $[0,T]$ with $T\le T^*$ as described by Theorem~1.1. Here $T^*$ is the maximal time of existence which is defined in the following sense:

\medskip

\noindent{\bf Definition} \em We call $T^*\in(0, \infty)$ to be the maximal time of existence of a smooth solution $(\rho-\tilde\rho,u,B)$ to \eqref{1.1}-\eqref{1.4} if for any $0<T<T^*$, $(\rho-\tilde\rho,u,B)$ solves \eqref{1.1}-\eqref{1.4} in $[0, T]\times\R^3$ and satisfies \eqref{1.1a}-\eqref{1.1b}; moreover, the conditions \eqref{1.1a}-\eqref{1.1b} fail to hold when $T=T^*$.
\rm 

\medskip

We will prove Theorem~1.2 using a contradiction argument. Therefore, for the sake of contradiction, we assume that

\begin{align}\label{2.1}
||\rho||_{L^\infty((0,T^*)\times\R^3)}+||\rho^{-1}||_{L^\infty((0,T^*)\times\R^3)}+||B||_{L^\infty((0,T^*)\times\R^3)}\le C.
\end{align}
\medskip

To facilitate our exposition, we first define some auxiliary functionals for $0\le t\le T\le T^*$:

\begin{align*}
A_1(t)=\sup_{0\le s \le t}\int_{\R^3}(|\nabla u|^2+|\nabla B|^2) dx+\int_0^t\!\!\!\int_{\R^3}(|\dot u|^2+|B_t|^2)dxds,
\end{align*}
\begin{align*}
A_2(t)=\sup_{0\le s \le t}\int_{\R^3}(|\dot u|^2+|B_t|^2) dx+\int_0^t\!\!\!\int_{\R^3}(|\nabla\dot u|^2+|\nabla B_t|^2)dxds,
\end{align*}
\begin{align*}
H(t)=\int_0^t\!\!\!\int_{\R^3}|\nabla u|^4dxds.
\end{align*}
\medskip

The following is the main theorem of this section:
\medskip

\noindent{\bf Theorem 2.1} \em Assume that the hypotheses and notations in  {\rm Theorem~1.1} are in force. Given $C>0$ and $\tilde\rho>0$, assume further that $(\rho-\tilde\rho,u,B)$ satisfies \eqref{2.1}. Then there exists a positive number $M$ which depends on $C_0,C,T^*$ and the system parameters $P,\varepsilon,\lambda,\nu$ such that, for $0\le t\le T\le T^*$,
\begin{align}\label{2.2}
A_1(t)+A_2(t)\le M.
\end{align}
\rm
\medskip

We prove Theorem~2.1 in a sequence of lemmas. We first derive the following lemma which gives estimates on the solutions of the Lam\'{e} operator $L=\varepsilon\Delta+(\varepsilon+\lambda)\nabla\divv$. More detailed discussions can also be found in Sun-Wang-Zhang \cite{sun}.

\medskip

\noindent{\bf Lemma 2.2} \em Consider the following equation:
\begin{equation}\label{2.3}
\varepsilon\Delta v+(\varepsilon+\lambda)\nabla\divv(v)=J,
\end{equation}
where $v=(v^1,v^2,v^3)(x)$, $J=(J^1,J^2,J^3)(x)$ with $x\in\R^3$ and $\varepsilon,\lambda>0$. Then for $p\in(1,\infty)$, we have:
\begin{enumerate}[label={\upshape{(2.\arabic*)}}, ref={\upshape{2.}\arabic*}, topsep=0.3cm, itemsep=0.15cm] 
\setcounter{enumi}{3}
\item if $J\in W^{2,p}(\R^3)$, then $||D^2_{x}v||_{L^p}\le \tilde C||J||_{L^p}$;
\item if $J=\nabla\phi$ with $\phi\in W^{2,p}(\R^3)$, then $||\nabla v||_{L^p}\le \tilde C||\phi||_{L^p}$;
\item if $J=\nabla\divv(\phi)$ with $\phi\in W^{2,p}(\R^3)$, then $||v||_{L^p}\le \tilde C||\phi||_{L^p}$.
\end{enumerate}
Here $\tilde C$ is a positive constant which depends only on $\varepsilon,\lambda,p$
\rm
\begin{proof}
A proof can be found in \cite{sun} pg. 39 and we omit the details here.
\end{proof}

We proceed to the following {\em a prior} estimates which is the energy-balanced law:
\medskip

\noindent{\bf Lemma 2.3} \em Assume that the hypotheses and notations of  {\rm Theorem~2.1} are in force. Then for any $0\le t\le T\le T^*$,
\begin{align}\label{2.7}
\sup_{0\le s \le t}\int_{\R^3}(|\rho-\tilde\rho|^2+\rho|u|^2+|B|^2) dx+\int_0^t\!\!\!\int_{\R^3}(|\nabla u|^2+|\nabla B|^2)dxds\le M(C)C_0,
\end{align}
where $M(C)$ is a constant which depends on $C$.
\rm
\begin{proof}
Let $G=G(\rho)$ be a functional defined by
\begin{align*}
G(\rho)=\rho\int_{\tilde\rho}^{\rho}s^{-1}(P(s)-P(\tilde\rho))ds.
\end{align*}
Multiplying the momentum equation \eqref{1.2} by $u^j$, summing over $j$, integrating and making use of the continuity equation \eqref{1.1}, we get:
\begin{align}\label{2.8}
\left.\int_{\R^3}\left[{\textstyle\frac{1}{2}}\rho|u|^{2}+G\right]dx\right |_{0}^{t} + \int_{0}^{t}\!\!\!&\int_{\R^3} u\cdot{\rm div}\,\left[({\textstyle\frac{1}{2}}|B|^{2})I_{3\times3}-BB^{T}\right]dxds\notag
\\ &+ \int_{0}^{t}\!\!\!\int_{\R^3}\left[\varepsilon|\nabla u|^{2} + (\varepsilon+\lambda)(\mathrm{div}\,u)^{2}\right]dxds=0.
\end{align}
Similarly, we multiply the magnetic field equation \eqref{1.3} by $B$ and integrate to get
\begin{equation}\label{2.9}
\left.\int_{\R^3}{\textstyle\frac{1}{2}}|B|^2 dx\right |_{0}^t + \int_{0}^t\!\!\!\int_{\R^3}B\cdot{\rm div}\left(Bu^T - uB^T\right)dxds=-\nu\int_{0}^t\!\!\!\int_{\R^3}|\nabla B|^2dxds.
\end{equation}
We then obtain \eqref{2.7} by adding \eqref{2.8} to \eqref{2.9} and using the fact that 
\begin{align*}
\int_{0}^t\!\!\!\int_{\R^3}\left[u\cdot{\rm div}[({\textstyle\frac{1}{2}}|B|^{2})I_{3\times3}-BB^{T}] + B\cdot{\rm div}(Bu^T - uB^T)\right]dxds=0.
\end{align*}
\end{proof}

We obtain the following $L^4$ bounds for $u$ and $B$:

\medskip

\noindent{\bf Lemma 2.4} \em Assume that the hypotheses and notations of  {\rm Theorem~2.1} are in force. Then for any $0\le t\le T\le T^*$,
\begin{align}\label{2.10}
\int_{\R^3}(|u(x,t)|^4+|B(x,t)|^4) dx\le M.
\end{align}
\rm
\begin{proof}
Multiply \eqref{1.2} by $2|u|^2 u$ and integrate to obtain
\begin{align}\label{2.11}
&\frac{d}{dt}\int_{\R^3}\rho|u|^4dx+\int_{\R^3}2|u|^2\big[\varepsilon|\nabla u|^2+(\lambda+\varepsilon)(\divv(u))^2\big]dx\notag\\
&\qquad\qquad\qquad\qquad+\int_{\R^3}8\left[\varepsilon|u|^{2}|\nabla(|u|)|^2+(\varepsilon+\lambda)(\divv(u))|u|u\cdot\nabla(|u|)\right] dx\notag\\
&=4\int_{\R^3}(P(\rho)-P(\tilde\rho))\divv(|u|^2 u)dx\notag\\
&\qquad\qquad\qquad\qquad+\int_{\R^3}2|B|^2\divv(|u|^2 u)dx+\int_{\R^3}4|u|^2u\cdot\divv(BB^T)dx
\end{align}
The third term on the left side of \eqref{2.11} can be estimated from below by
\begin{align*}
&\int_{\R^3}8[\varepsilon|u|^{2}|\nabla(|u|)|^2+(\varepsilon+\lambda)(\divv(u))|u|u\cdot\nabla(|u|)] dx\\
&\qquad\qquad\qquad\qquad\qquad\qquad\ge\int_{\R^3}4|u|^2\left[\varepsilon|\nabla u|^2+2(\varepsilon-\frac{\varepsilon+\lambda}{2})|\nabla(|u|)|^2\right]dx.
\end{align*}
By assumption \eqref{1.7} we have $\varepsilon<\lambda$, hence it implies
\begin{align}\label{2.12}
\int_{\R^3}8|u|^2(\varepsilon-\frac{\varepsilon+\lambda}{2})|\nabla(|u|)|^2\ge M\int_{\R^3}|u|^2|\nabla u|^2dx.
\end{align}
On the other hand, we multiply \eqref{1.3} by $4|B|^2 B$ and integrate to get
\begin{align}\label{2.13}
&\frac{d}{dt}\int_{\R^3}|B|^2 dx+\int_{\R^3}4\nu|B|^2|\nabla B|^2 dx+\int_{\R^3}2\nu |B|^2|\nabla B|^2dx\notag\\
&\qquad\qquad\qquad\qquad\qquad\qquad=-\int_{\R^3}|B|^2 B\cdot\divv(Bu^T-uB^T)dx.
\end{align}
Adding \eqref{2.13} to \eqref{2.12} and integrate with respect to $t$, we get
\begin{align}\label{2.14}
&\left(\int_{\R^3}(\rho|u|^4+|B|^4)dx\right)+\int_0^t\!\!\!\int_{\R^3}(|u|^2|\nabla u|^2+|B|^2|\nabla B|^2)dxds\notag\\
&\le M\left[\int_0^t\!\!\!\int_{\R^3}4(P(\rho)-P(\tilde\rho))\divv(|u|^2u)dxds+\int_0^t\!\!\!\int_{\R^3}2|B|^2\divv(|u|^2u)dxds\right]\notag\\
&\qquad-M\left[\int_0^t\!\!\!\int_{\R^3}2|u|^2u\cdot\divv(BB^T)dxds+\int_0^t\!\!\!\int_{\R^3}|B|^2B\cdot\divv(Bu^T-uB^T)dxds\right].
\end{align}
Using the assumption \eqref{2.1}, the right side of \eqref{2.14} can be bounded by
\begin{align}\label{2.15}
&\left[\int_0^t\!\!\!\int_{\R^3}(|P(\rho)-P(\tilde\rho)|^2|u|^2+|B|^4|u|^2)dxds\right]^{\frac{1}{2}}\left[\int_0^t\!\!\!\int_{\R^3}(|u|^2|\nabla u|^2+|B|^2|\nabla B|^2)dxds\right]^{\frac{1}{2}}\notag\\
&\qquad\le M\left[T^*C_0+\int_0^t\!\!\!\int_{\R^3}(|B|^4+|u|^4)dxds\right]^{\frac{1}{2}}\left[\int_0^t\!\!\!\int_{\R^3}(|u|^2|\nabla u|^2+|B|^2|\nabla B|^2)dxds\right]^{\frac{1}{2}}.
\end{align}
Using \eqref{2.15} on \eqref{2.14} and applying Cauchy Inequality, we get
\begin{align*}
\int_{\R^3}(|u|^4+|B|^4)dx\le M+\int_0^t\!\!\!\int_{\R^3}(|u|^4+|B|^4)dxds,
\end{align*}
and \eqref{2.10} now follows by Gronwall's inequality.
\end{proof}

We obtain estimates on the functional $A_1$ in terms of $H$:

\medskip

\noindent{\bf Lemma 2.5} \em Assume that the hypotheses and notations of  {\rm Theorem~2.1} are in force. Then for any $0\le t\le T\le T^*$,
\begin{align}\label{2.16}
A_1(t)\le M[1+H(t)].
\end{align}
\rm
\begin{proof}
We multiply \eqref{1.2} by $\dot{u}^j$, sum over $j$ and integrate to get
\begin{align}\label{2.17}
&\int_{\R^3}|\nabla u|^2dx+\int_0^t\!\!\!\int_{\R^3}\rho|\dot u|^2dxds\notag\\
&\qquad\le C_0+\left|\int_0^t\!\!\!\int_{\R^3}[\dot u\cdot\nabla(\frac{1}{2}|B|^2)-\dot{u}\cdot\divv(BB^T)]\right|+\int_0^t\!\!\!\int_{\R^3}|\nabla u|^3dxds.
\end{align}
Niext we multiply \eqref{1.3} by $B_t$ and integrate,
\begin{align}\label{2.18}
\int_{\R^3}|\nabla B|^2dx+\int_0^t\!\!\!\int_{\R^3}|B_t|^2dxds
\le C_0+\left|\int_0^t\!\!\!\int_{\R^3}B_t \cdot\divv(uB^T-u^TB)dxds\right|.
\end{align}
Adding \eqref{2.17} and \eqref{2.18}, we obtain
\begin{align}\label{2.19}
&\int_{\R^3}(|\nabla u|^2+|\nabla B|^2)dx+\int_0^t\!\!\!\int_{\R^3}(\rho|\dot u|^2+|B_t|^2)dxds\notag\\
&\le C_0+\int_0^t\!\!\!\int_{\R^3}|\nabla u|^3dxds+\int_0^t\!\!\!\int_{\R^3}(|\nabla B|^2|B|^2+|\nabla u|^2|B|^2+|\nabla B|^2|u|^2)dxds.
\end{align}
The second term on the right side of \eqref{2.19} is bounded by
\begin{align*}
\left(\int_0^t\!\!\!\int_{\R^3}|\nabla u|^2dxds\right)^\frac{1}{2}\left(\int_0^t\!\!\!\int_{\R^3}|\nabla u|^4dxds\right)^\frac{1}{2}\le C_0+H(t),
\end{align*}
where the last inequality follows by Lemma~2.3. For the last integral on the right side of \eqref{2.19}, using assumption \eqref{2.1}, it can be bounded by $\dis\int_0^t\!\!\!\int_{\R^3}|\nabla B|^2|u|^2dxds+M\int_0^t\!\!\!\int_{\R^3}(|\nabla u|^2+|\nabla B|^2)dxds$. So it remains to estimate $\dis\int_0^t\!\!\!\int_{\R^3}|\nabla B|^2|u|^2dxds$.
\medskip

Recall from Lemma~2.4 that, for $0\le t\le T\le T^*$,
\begin{align*}
\int_{\R^3}|u|^4dx\le M,
\end{align*}
Therefore, using \eqref{1.9},
\begin{align*}
\int_0^t\!\!\!\int_{\R^3}|\nabla B|^2|u|^2dxds&\le \int_0^t\left(\int_{\R^3}|\nabla B|^4dx\right)^\frac{1}{2}\left(\int_{\R^3}|u|^4dx\right)^\frac{1}{2}ds\\
&\le M\int_0^t\left(\int_{\R^3}|\nabla B|^2dx\right)^\frac{1}{4}\left(\int_{\R^3}|D^2_{x}B|^2dx\right)^\frac{3}{4}ds\\
&\le M\left(\int_0^t\!\!\!\int_{\R^3}|\nabla B|^2dxds\right)^\frac{1}{4}\left(\int_0^t\!\!\!\int_{\R^3}|\Delta B|^2dxds\right)^\frac{3}{4}\\
&\le MC_0^\frac{1}{4}\left[\int_0^t\!\!\!\int_{\R^3}(|B_t|^2+|\nabla B|^2|u|^2+|\nabla u|^2|B|^2)dxds\right]^\frac{3}{4}\\
&\le MC_0^\frac{1}{4}\left[A_1(t)+C_0+\int_0^t\!\!\!\int_{\R^3}|\nabla B|^2|u|^2dxds\right]^\frac{3}{4},
\end{align*}
and by Cauchy inequality, we obtain
\begin{align}\label{2.20}
\int_0^t\!\!\!\int_{\R^3}|\nabla B|^2|u|^2dxds\le M\left(1+A_1(t)^\frac{3}{4}\right).
\end{align}
Applying \eqref{2.20} to \eqref{2.19} and absorbing terms, \eqref{2.16} follows.
\end{proof}

We derive the following estimates on the {\em effective viscous flux} which were first described by Hoff \cite{hoff1} and later modified by Sun-Wang-Zhang \cite{sun}.

\medskip

\noindent{\bf Lemma 2.6} \em Assume that the hypotheses and notations of  {\rm Theorem~2.1} are in force. Then for any $0\le t\le T\le T^*$,
\begin{align}\label{2.21}
\sup_{0\le s\le t}\int_{\R^3}|\nabla w|^2+\int_0^t\!\!\!\int_{\R^3}|w_t|^2dxds+\int_0^t\!\!\!\int_{\R^3}|D_{x}^2 w|^2dxds\le M,
\end{align}
where $w=u-v$ with v satisfying:
\begin{align*}
\varepsilon\Delta v+(\varepsilon+\lambda)\nabla\divv(v)=\nabla(P(\rho)-P(\tilde\rho)).
\end{align*}
\rm
\begin{proof}
Using the momentum eqaution \eqref{1.2},
\begin{align}\label{2.22}
\rho w_t-\varepsilon\Delta w-(\varepsilon+\lambda)\nabla\divv(w)=-\rho u\cdot\nabla u-\rho v_t-\nabla(\frac{1}{2}|B|^2)+\divv(BB^T).
\end{align}
Multiply \eqref{2.22} by $w_t$ and integrate,
\begin{align}\label{2.23}
&\left.\int_{\R^3}\varepsilon|\nabla w|^2dx\right|_0^t+\int_0^t\!\!\!\int_{\R^3}(\varepsilon+\lambda)|\divv w|^2dxds+\int_0^t\!\!\!\int_{\R^3}\rho|w_t|^2dxds\notag\\
&\qquad\qquad\qquad=\int_0^t\!\!\!\int_{\R^3}[-\rho u\cdot\nabla u-\rho v_t-\nabla(\frac{1}{2}|B|^2)+\divv(BB^T)]\cdot w_t dxds.
\end{align}
The first term on the right side of \eqref{2.23} can be estimated as follows:
\begin{align}\label{2.24}
\int_0^t\!\!\!\int_{\R^3}(-\rho u\cdot\nabla u)\cdot w_t dxds\le M\left[\int_0^t\!\!\!\int_{\R^3}|u|^2|\nabla u|^2 dxds\right]^\frac{1}{2}\left[\int_0^t\!\!\!\int_{\R^3}|w_t|^2 dxds\right]^\frac{1}{2}.
\end{align}
For the term $\dis\int_0^t\!\!\!\int_{\R^3}|u|^2|\nabla u|^2 dxds$, using Lemma~2.4,
\begin{align}\label{2.24a}
\int_0^t\!\!\!\int_{\R^3}|u|^2|\nabla u|^2 dxds&\le\int_0^t\left(\int_{\R^3}|u|^4dx\right)^\frac{1}{2}\left(\int_{\R^3}|\nabla u|^4dx\right)^\frac{1}{2}\notag\\
&\le M\int_0^t\left(\int_{\R^3}(|\nabla w|^4+|\nabla v|^4)dx\right)^\frac{1}{2}ds\notag\\
&\le M\left[\int_0^t\left(\int_{\R^3}|\nabla w|^2dx\right)^\frac{1}{4}\left(\int_{\R^3}|D_{x}^2 w|^2dx\right)^\frac{3}{4}ds\right]\notag\\
&\qquad\qquad\qquad\qquad\qquad\qquad+M\int_0^t\left(\int_{\R^3}|\nabla v|^4 dx\right)^\frac{1}{2}\notag\\
&\le M\left(\int_0^t\!\!\!\int_{\R^3}|\nabla w|^2dxds\right)^\frac{1}{4}\left(\int_0^t\!\!\!\int_{\R^3}|D_{x}^2 w|^2dxds\right)^\frac{3}{4}\notag\\
&\qquad\qquad\qquad\qquad+M\int_0^t\left(\int_{\R^3}|P(\rho)-P(\tilde\rho)|^4dx\right)^\frac{1}{2}ds\notag\\
&\le M\left(\int_0^t\!\!\!\int_{\R^3}|\nabla w|^2dxds\right)^\frac{1}{4}\left(\int_0^t\!\!\!\int_{\R^3}|D_{x}^2 w|^2dxds\right)^\frac{3}{4}\notag\\
&\qquad\qquad\qquad\qquad\qquad+MC_0^\frac{1}{2}{T^*}^\frac{1}{2},
\end{align}
where the last inequality follows by Lemma~2.2 and assumption \eqref{2.1}. Therefore \eqref{2.24} becomes
\begin{align}\label{2.25}
&\int_0^t\!\!\!\int_{\R^3}(-\rho u\cdot\nabla u)\cdot w_t dxds\notag\\
&\le M\left[\left(\int_0^t\!\!\!\int_{\R^3}|\nabla w|^2dxds\right)^\frac{1}{4}\left(\int_0^t\!\!\!\int_{\R^3}|D_{x}^2 w|^2dxds\right)^\frac{3}{4}+MC_0^\frac{1}{2}{T^*}^\frac{1}{2}\right]^\frac{1}{2}\notag\\
&\qquad\qquad\qquad\qquad\qquad\qquad\qquad\qquad\qquad\qquad\qquad\times\left[\int_0^t\!\!\!\int_{\R^3}|w_t|^2 dxds\right]^\frac{1}{2}.
\end{align}
The third and the fourth term on the right side of \eqref{2.23} are bounded by 
\begin{align}\label{2.26}
\int_0^t\!\!\!\int_{\R^3}|\nabla B||B||w_t|dxds&\le M\left(\int_0^t\!\!\!\int_{\R^3}|\nabla B|^2|B|^2dxds\right)^\frac{1}{2}\left(\int_0^t\!\!\!\int_{\R^3}|w_t|^2dxds\right)^\frac{1}{2}\notag\\
&\le MC_0^\frac{1}{2}\left(\int_0^t\!\!\!\int_{\R^3}|w_t|^2dxds\right)^\frac{1}{2}.
\end{align}
It remains to estimate the term $\dis\int_0^t\!\!\!\int_{\R^3}-\rho v_t\cdot w_t dxds$ on the right side of \eqref{2.23}. By the definition of $v$ and $P(\rho)$, we have
\begin{align*}
\varepsilon\Delta v_t+(\varepsilon+\lambda)\nabla\divv(v_t)&=\nabla P(\rho)_t=\nabla\divv(-P(\rho)u),
\end{align*}
Hence we can apply Lemma~2.2 and Lemma~2.3 to get
\begin{align*}
\int_0^t\!\!\!\int_{\R^3}|v_t|^2 dxds\le \int_0^t\!\!\!\int_{\R^3}|P(\rho)u|^2dxds\le MC_0, 
\end{align*}
and 
\begin{align}\label{2.27}
\int_0^t\!\!\!\int_{\R^3}-\rho v_t\cdot w_t dxds&\le M\left(\int_0^t\!\!\!\int_{\R^3}|w_t|^2 dxds\right)^\frac{1}{2}\left(\int_0^t\!\!\!\int_{\R^3}|v_t|^2 dxds\right)^\frac{1}{2}\notag\\
&\le M\left(\int_0^t\!\!\!\int_{\R^3}|w_t|^2 dxds\right)^\frac{1}{2}.
\end{align}
Using \eqref{2.25}, \eqref{2.26} and \eqref{2.27} on \eqref{2.23},
\begin{align}\label{2.28}
&\left.\int_{\R^3}\varepsilon|\nabla w|^2dx\right|_0^t+\int_0^t\!\!\!\int_{\R^3}(\varepsilon+\lambda)|\divv(w)|^2dxds+\int_0^t\!\!\!\int_{\R^3}\rho|w_t|^2 dxds\notag\\
&\qquad\qquad\qquad\le M\left(\int_0^t\!\!\!\int_{\R^3}|\nabla w|^2 dxds\right)^\frac{1}{4}\left(\int_0^t\!\!\!\int_{\R^3}|D_{x}^2 w|^2 dxds\right)^\frac{3}{4}+M.
\end{align}
It remains to estimate $\dis\int_0^t\!\!\!\int_{\R^3}|D_{x}^2 w|^2 dxds$. We rearrange the terms in \eqref{2.22} to get
\begin{align*}
\varepsilon\Delta w+(\varepsilon+\lambda)\nabla\divv(w)=\rho w_t+\rho \nabla u\cdot u+\rho v_t+\nabla(\frac{1}{2}|B|^2)-\divv(BB^T),
\end{align*}
and so by Lemma~2.2,
\begin{align*}
\int_0^t\!\!\!\int_{\R^3}|D_{x}^2 w|^2 dxds\le\int_0^t\!\!\!\int_{\R^3}(|\rho w_t|^2+|\rho\nabla u\cdot u|^2+|\rho v_t|^2+|\nabla B|^2|B|^2) dxds\notag\\
\le M\left[\int_0^t\!\!\!\int_{\R^3}|w_t|^2 dxds+\left(\int_0^t\!\!\!\int_{\R^3}|\nabla w|^2 dxds\right)^\frac{1}{4}\left(\int_0^t\!\!\!\int_{\R^3}|D_{x}^2 w|^2 dxds\right)^\frac{3}{4}+1\right].
\end{align*}
Therefore
\begin{align*}
\int_0^t\!\!\!\int_{\R^3}|D_{x}^2 w|^2 dxds\le M\left[\int_0^t\!\!\!\int_{\R^3}|w_t|^2 dxds+\int_0^t\!\!\!\int_{\R^3}|\nabla w|^2 dxds+1\right],
\end{align*}
and we apply the above to \eqref{2.28} to conclude
\begin{align*}
\left.\int_{\R^3}\varepsilon|\nabla w|^2 dx\right|_0^t+\int_0^t\!\!\!\int_{\R^3}(\varepsilon+\lambda)|\divv(w)|^2dxds+\int_0^t\!\!\!\int_{\R^3}\rho|w_t|^2 dxds\\
\le M\left[\int_0^t\!\!\!\int_{\R^3}|\nabla w|^2dxds+1\right],
\end{align*}
which \eqref{2.21} follows by Gronwall's inequality.
\end{proof}

We finally obtain an estimate on the functional $A_2$ which is sufficient to prove Theorem~2.1:

\medskip

\noindent{\bf Lemma 2.7} \em Assume that the hypotheses and notations of  {\rm Theorem~2.1} are in force. Then for any $0\le t\le T\le T^*$,
\begin{align}\label{2.29}
A_2(t)\le M\left[A_1(t)+H(t)+1\right]
\end{align}
\rm
\begin{proof}
Taking the convective derivative in the momentum equation \eqref{1.2}, multiplying it by $\dot{u}^j$, summing over $j$ and integrating,
\begin{align}\label{2.30}
&\sup_{0\le s\le t}\int_{\R^3}|\dot{u}|^2dx+\int_0^t\!\!\!\int_{\R^3}|\nabla\dot u|^2dxds\notag\\
&\qquad\qquad\qquad\le C_0+H(t)+\int_0^t\!\!\!\int_{\R^3}|B|^2(|B_t|^2+|u|^2|\nabla B|^2)dxds.
\end{align}
Next we differentiate the magnetic field equation \eqref{1.3} with respect to $t$, multiply by $B_t$ and integrate,
\begin{align*}
&\left.\frac{1}{2}\int_{\R^3}|B_t|^2dx\right|_0^t+ \nu\int_{0}^{t}\!\!\!\int_{\R^3}|\nabla B_t|^2 dxds\\
&\qquad\qquad\qquad=-\int_{0}^{t}\!\!\!\int_{\R^3} B_t\cdot[\divv(B u^{T}-u B^{T})]_{t} dxds.
\end{align*}
Adding the above to \eqref{2.30} and absorbing terms,
\begin{align}\label{2.31}
&\sup_{0\le s\le t}\int_{\R^3}(|\dot u|^2+|B_t|^2)dx+\int_0^t\!\!\!\int_{\R^3}(|\nabla\dot u|^2+|\nabla B_t|^2) dxds\notag\\
&\qquad\qquad\le M\left[A_1+H+\int_0^t\!\!\!\int_{\R^3}|B|^2|u|^2(|\nabla u|^2+|\nabla B|^2)dxds\right]\notag\\
&\qquad\qquad\qquad\qquad+M\int_0^t\!\!\!\int_{\R^3}(|B|^2|B_t|^2+|B|^2|\dot u|^2+|B_t|^2|u|^2)dxds
\end{align}
The third term on the right side of \eqref{2.31} is bounded by
\begin{align}\label{2.32}
\int_0^t\!\!\!\int_{\R^3}|u|^2(|\nabla u|^2+|\nabla B|^2)dxds&\le \left[\int_0^t\!\!\!\int_{\R^3}|u|^2|\nabla u|^2dxds+\int_0^t\!\!\!\int_{\R^3}|u|^2|\nabla B|^2dxds\right]\notag\\
&\le M\left[1+A_1^{\frac{3}{4}}\right],
\end{align}
where the last inequality follows by \eqref{2.20} and \eqref{2.24a}. The last term on the right side of \eqref{2.31} is bounded by
\begin{align}\label{2.33}
\int_0^t\!\!\!\int_{\R^3}(|B_|^2+|\dot u|^2)dxds&+\int_0^t\!\!\!\int_{\R^3}|B_t|^2|u|^2dxds\notag\\
&\le MA_1+\int_0^t\left(\int_{\R^3}|B_t|^4dx\right)^\frac{1}{2}\left(\int_{\R^3}|u|^4dx\right)^\frac{1}{2}ds\notag\\
&\le MA_1+M\int_0^t\left(\int_{\R^3}|B_t|^2dx\right)^\frac{1}{4}\left(\int_{\R^3}|\nabla B_t|^2dx\right)^\frac{3}{4}ds\notag\\
&\le MA_1+M\left(\int_0^t\!\!\!\int_{\R^3}|B_t|^2dxds\right)^\frac{1}{4}\left(\int_0^t\!\!\!\int_{\R^3}|\nabla B_t|^2dxds\right)^\frac{3}{4}\notag\\
&\le MA_1+A_1^\frac{1}{4}A_2^\frac{3}{4}.
\end{align}
Using \eqref{2.32} and \eqref{2.33} on \eqref{2.31} and absorbing terms, \eqref{2.30} follows.
\end{proof}
\begin{proof}[proof of Theorem~2.1]
Recall from Lemma~2.5 and 2.7 that
\begin{align*}
A_1\le M\left[H+1\right]
\end{align*}
and
\begin{align*}
A_2\le M\left[H+A_1\right].
\end{align*}
So it remains to estimate $H$. Let $w$ and $v$ be as defined in Lemma~2.6. Then
\begin{align}\label{2.34}
\int_0^t\!\!\!\int_{\R^3}|\nabla u|^4dxds\le \int_0^t\!\!\!\int_{\R^3}|\nabla w|^4dxds+\int_0^t\!\!\!\int_{\R^3}|\nabla v|^4 dxds
\end{align}
The second term on the right side of \eqref{2.34} is bounded by $\dis\int_0^t\!\!\!\int_{\R^3}|P(\rho)-P(\tilde\rho)|^4dxds$. And for $\dis\int_0^t\!\!\!\int_{\R^3}|\nabla w|^4dxds$, using \eqref{2.21},
\begin{align*}
\int_0^t\!\!\!\int_{\R^3}|\nabla w|^4dxds&\le\int_0^t\left(\int_{\R^3}|\nabla w|^2dx\right)^\frac{1}{2}\left(\int_{\R^3}|D^2_{x} w|^2dx\right)^\frac{3}{2}\\
&\le\left(\sup_{0\le s\le t}\int_{\R^3}|D^2_{x} w|^2dx\right)^\frac{1}{2}\left(\sup_{0\le s\le t}\int_{\R^3}|\nabla w|^2dx\right)^\frac{1}{2}\left(\int_0^t\!\!\!\int_{\R^3}|\nabla w|^2dx\right)\\
&\le M\left(\sup_{0\le s\le t}\int_{\R^3}|D^2_{x} w|^2dx\right)^\frac{1}{2}.
\end{align*}
Notice that, by rearranging the terms in \eqref{2.22},
\begin{align*}
\varepsilon\Delta w+(\varepsilon+\lambda)\nabla\divv(w)=\rho\dot{u}+\nabla(\frac{1}{2}|B|^2)-\divv(BB^T),
\end{align*}
and so by Lemma~2.2,
\begin{align*}
\int_{\R^3}|D^2_{x} w|^2dx&\le\int_{\R^3}|\rho\dot u|^2dx+\int_{\R^3}|\nabla(\frac{1}{2}|B|^2)-\divv(BB^T)|^2dx\\
&\le M\left[\int_{\R^3}|\dot u|^2dx+\int_{\R^3}|\nabla B|^2|B|^2dx\right]\\
&\le M[A_2+A_1].
\end{align*}
Therefore we conclude that
\begin{align*}
H\le M\left[A_2+A_1\right]^\frac{1}{2},
\end{align*}
and \eqref{2.2} follows.
\end{proof}

\bigskip

\section{Higher Order Estimates and proof of Theorem~1.2}

\bigskip

In this section we continue to obtain higher order estimates on the smooth local solution $(\rho-\tilde\rho,u,B)$ as described in section~2. Together with Theorem~2.1, we show that, under the assumption \eqref{2.1}, the smooth local solution to \eqref{1.1}-\eqref{1.4} can be extended beyond the maximal time of existence $T^*$ as defined in section~2, thereby contradicting the maximality of $T^*$. The following is the main theorem of this section:

\medskip

\noindent{\bf Theorem 3.1} \em Assume that the hypotheses and notations in  {\rm Theorem~2.1} are in force. Given $C>0$ and $\tilde\rho>0$, assume further that $(\rho-\tilde\rho,u,B)$ satisfies \eqref{2.1}. Then there exists a positive number $M'$ which depends on $C_0,C,T^*$ and the system parameters $P,\varepsilon,\lambda,\nu$ such that, for $0\le t\le T\le T^*$,
\begin{align}\label{3.1}
\sup_{0\le s\le t}||(\rho-\tilde\rho,u,B)||_{H^3(\R^3)}+\int_0^t||(u,B)(\cdot,s)||^2_{H^4(\R^3)}ds\le M'
\end{align}
\rm
\medskip

\noindent{\em Proof. }
We give the proof in a sequence of steps. Most of the details are reminiscent of Suen and Hoff \cite{suenhoff} and we omit those which are identical to or nearly identical to arguments given in \cite{suenhoff}. We first begin with the following estimates on the effective viscous flux $F$ and the vorticity matrix $\omega$:

\medskip

\noindent{\bf  Step 1:} {\em Define}
\begin{align*}
F&=(2\varepsilon+\lambda)\divv(u)-(P(\rho)-P(\tilde\rho)),\\
\omega&=\omega^{j,k}=u^j_{x_k}-u^k_{x_j},
\end{align*}
{\em Then for} $q\in(1,\infty)$,
\begin{align}\label{3.2}
||\nabla u(\cdot,t)||_{L^q}&\le M(q)\left[||F(\cdot,t)||_{L^q}+||\omega(\cdot,t)||_{L^q}+||(P(\rho)-P(\tilde\rho))(\cdot,t)||_{L^q}\right],\\
||\nabla\omega(\cdot,t)||_{L^q}&\le (q)M[||\rho\dot u(\cdot,t)||_{L^q}+||\nabla B\cdot B(\cdot,t)||_{L^q}],\label{3.3}
\end{align}
{\em where $M(q)$ is a positive constant depending on $q$ and}
\begin{align}\label{3.4}
\sup_{0\le s\le t}\int_{\R^3}(|\nabla F|^2+|\nabla\omega|^2)dx\le M'.
\end{align}
\begin{proof}[proof of Step 1]
We give the proof of \eqref{3.2} as an example. Using the definition of $F$ and $\omega$,
\begin{align}
(2\varepsilon+\lambda)\Delta u^j&=F_{x_j}+(2\varepsilon+\lambda)\omega^{j,k}_{x_k}+(P(\rho)-P(\tilde\rho)_{x_j},\label{3.5}
\end{align}
Differentiating and taking the Fourier transform we then obtain
\begin{align*}
(2\varepsilon+\lambda)\hat{u}^j_{x_l}(y,t)=\frac{y_j y_l}{|y|^2}\hat{F}(y,t)+(2\varepsilon+\lambda)\frac{y_k y_l}{|y|^2}\widehat{\omega^{j,k}}(y,t)+\frac{y_k y_l}{|y|^2}(\widehat{P-\tilde{P}})(y,t)
\end{align*}
and \eqref{3.2} then follows immediately from the Marcinkiewicz multiplier theorem (Stein \cite{stein}, pg. 96). Similarly, \eqref{3.3} can be proved by the same method. Also, by the definition of $F$, we have
\begin{align}
\Delta F=\divv(g),
\end{align} 
where $g^j=\rho\dot u^j+(\frac{1}{2}|B|^2)_{x_j}-\divv(B^j B)$. So we have
\begin{align*}
\sup_{0\le s\le t}\int_{\R^3}|\nabla F|^2dx&\le\sup_{0\le s\le t}\int_{\R^3}|g|^2dx\\
&\le\sup_{0\le s\le t}\int_{\R^3}(\rho|\dot u|^2+|\nabla B|^2|B|^2)dx\le M',
\end{align*}
and similarly, $\dis\sup_{0\le s\le t}\int_{\R^3}|\nabla\omega|^2dx\le M'$, which proves \eqref{3.4}.
\end{proof}
\noindent{\bf  Step 2:} {\em The velocity gradient satisfies the following bound}
\begin{align*}
\int_0^t||\nabla u(\cdot,t)||_{L^\infty}ds\le M'.
\end{align*}
\begin{proof}[proof of Step 2]
The proof is identical to Suen and Hoff \cite{suenhoff} pg. 51--53, and we omit the details here.
\end{proof}
\noindent{\bf  Step 3:} {\em We further obtain}
\begin{align}
||D^2_{x}u(\cdot,t)||_{L^2}&\le M'\left[||\rho u\dot(\cdot,t)||_{L^2}+||\nabla B\cdot B(\cdot,t)||_{L^2}+||\nabla P(\cdot,t)||_{L^2}\right],\label{3.5a}
\end{align}
\begin{align}
||D^3_{x}u(\cdot,t)||_{L^2}&\le M'\left[||\nabla\rho\cdot\dot u(\cdot,t)||_{L^2}+||\rho\nabla\dot u(\cdot,t)||_{L^2}+||B\cdot D^2_{x}B(\cdot,t)||_{L^2}\right]\notag\\
&\qquad+M'\left[|||\nabla B|^2(\cdot,t)||_{L^2}+||D^2_{x}P(\cdot,t)||_{L^2}\right].\label{3.5b}
\end{align}
\begin{proof}[proof of Step 3]
These follow immediately from the momentum equation \eqref{1.2} and the ellipticity of the Lam\'{e} operator $\varepsilon\Delta+(\varepsilon+\lambda)\nabla\divv$.
\end{proof}
\noindent{\bf  Step 4:} {\em The following $H^2$-bound for density holds}
\begin{align}\label{3.5c}
\sup_{0\le s\le t}||(\rho-\tilde\rho)(\cdot,s)||_{H^2}\le M'.
\end{align}
\begin{proof}[proof of Step 4]
We take the spatial gradient of the mass equation \eqref{1.1}, multiply by $\nabla\rho$ and integrate by parts to obtain
\begin{align}
\frac{\partial}{\partial t}\int_{\R^3}|\nabla\rho|^2dx\le M'\left[\int_{\R^3}|\nabla\rho|^2dx+\int_{\R^3}|D^2_{x} u|^2dx\right]\label{3.6}
\end{align}
From \eqref{3.5a},
\begin{align*}
\int_0^t\!\!\!\int_{\R^3}|D^2_{x} u|^2dxds&\le\int_0^t\!\!\!\int_{\R^3}(|\dot u|^2+|\nabla B\cdot B|^2+|\nabla\rho|^2)dxds\\
&\le M'+\int_0^t\!\!\!\int_{\R^3}|\nabla\rho|^2dxds.
\end{align*}
Applying the above to \eqref{3.6} and using the result of Step~2,
\begin{align*}
\sup_{0\le s\le t}||\nabla\rho(\cdot,s)||_{L^2}\le M'.
\end{align*}
By similar argument, we can show that $\dis\sup_{0\le s\le t}||D^2_{x}\rho(\cdot,s)||_{L^2}\le M'$ and \eqref{3.5c} follows.
\end{proof}
\noindent{\bf  Step 5:} {\em The velocity and magnetic field satisfy}
\begin{align}\label{3.7}
\sup_{0\le s\le t}\left(||u(\cdot,s)||_{H^3}+||B(\cdot,s)||_{H^3}\right)\le M'.
\end{align}
\begin{proof}[proof of Step 5]
Define the forward difference of quotient $D^h_{t}$ by
\begin{align*}
D^h_t(f)(t)=(f(t+h)-f(t))h^{-1}
\end{align*}
and let $E^j=D^h_t(u^j)+u\cdot\nabla u^j$. By differentiating the momentum equation, we obtain
\begin{align*}
&\int_{\R^3}\rho|E_{x_j}|^2dx+\int_0^t\!\!\!\int_{\R^3}\left(|\nabla E_{x_j}|^2+|D^h_t(\divv(u_{x_j})+u\cdot\nabla(\divv(u_{x_j})|^2\right)dxds\\
&\qquad\qquad\qquad\le M'+\int_0^t\!\!\!\int_{\R^3}|\nabla E|^2dxds+\mathcal O(h),
\end{align*}
where $\mathcal O(h)\rightarrow 0$ as $h\rightarrow0$. Therefore by taking $h\rightarrow0$,
\begin{align*}
\sup_{0\le s\le t}||\nabla\dot u(\cdot,s)||_{L^2}+\int_0^t\!\!\!\int_{\R^3}|D^2_x\dot u|^2dxds\le M'.
\end{align*}
The bound for $\nabla B_t$ can be derived in an exactly same way.
\end{proof}
\noindent{\bf  Step 6:} {\em Finally we have the following bounds}
\begin{align}
&\int_0^t\!\!\!\int_{\R^3}(|D^4_{x}u|^2+|D^4_{x}B|^2)dxds\le M'\left[1+\int_0^t\!\!\!\int_{\R^3}|D^3_{x}\rho|^2dxds\right],\label{3.8a}\\
&\sup_{0\le s\le t}\left(||D^3_{x}\rho(\cdot,s)||_{L^2}+||D^3_{x}B(\cdot,s)||_{L^2}\right)+\int_0^t\!\!\!\int_{\R^3}|D^4_{x}u|^2dxds\le M'.\label{3.8b}
\end{align}
\begin{proof}[proof of Step 6]
For \eqref{3.8a}, it can be obtained by differentiating \eqref{1.2} and \eqref{1.3} twice with respect to space, expressing the fourth derivatives of $u$ and $B$ in the terms second derivatives of $\dot u$, $B_t$, $\nabla\rho$ and lower order terms, and applying the bounds in \eqref{2.1} and \eqref{3.7}

For \eqref{3.8b}, it can be obtained by applying two space derivatives and one spatial difference operator $D^h_{x_j}$ defined by
\begin{align*}
D^h_{x_j}(f)(t)=(f(x+he_j)-f(x))h^{-1}
\end{align*}
such that
\begin{align*}
\int_{\R^3}|D^h_{x_j} D_{x_i} D_{x_k}\rho|^2dx&\le M'+\int_0^t\!\!\!\int_{\R^3}(|D^4_x u|^2+|D^h_{x_j}D_{x_i}D_{x_k}\rho|^2)dxds\\
&\le M'+\int_0^t\!\!\!\int_{\R^3}|D^3_{x}\rho|^2dxds.
\end{align*}
Taking $h\rightarrow0$ and applying Gronwall's inequality, we obtain the required bound for the term $\dis||D^3_{x}\rho(\cdot,s)||_{L^2}$. 
\end{proof}

\begin{proof}[proof of Theorem~1.2]
Using Theorem~3.1, we can apply an open-closed argument on the time interval which is identical to the one given in Suen and Hoff \cite{suenhoff} pp. 31 to extend the local solution $(\rho-\tilde\rho,u,B)$ beyond $T^*$, which contradicts the maximality of $T^*$. Therefore the assumption \eqref{2.1} does not hold and this completes the proof of Theorem~1.2.
\end{proof}



\end{document}